\documentclass[11pt]{amsart}
\usepackage[T1]{fontenc}
\usepackage{lmodern}
\usepackage[letterpaper,margin=1.06in]{geometry}
\usepackage{amsmath,amssymb,amsthm,mathtools,mathrsfs}
\usepackage{microtype}
\usepackage{graphicx}
\usepackage[hidelinks]{hyperref}
\usepackage{enumitem}
\usepackage{tikz}
\usetikzlibrary{arrows.meta,positioning}
\setlist[enumerate]{itemsep=2pt,topsep=3pt}
\allowdisplaybreaks[2]
\DeclareMathOperator{\Tr}{Tr}
\DeclareMathOperator{\Ran}{Ran}
\DeclareMathOperator{\diag}{diag}
\DeclareMathOperator{\sgn}{sgn}

\newcommand{\C}{\mathbb C}
\newcommand{\R}{\mathbb R}
\newcommand{\In}{\operatorname{In}}
\newcommand{\supp}{\operatorname{supp}}

\newcommand{\RR}{\mathsf{R}}
\newcommand{\XX}{\mathsf X}
\newcommand{\EE}{\mathsf E}
\newcommand{\GK}{\mathscr G}
\newcommand{\TT}{\mathscr T}
\newcommand{\BB}{\mathscr B}
\newcommand{\CF}{\mathscr C}
\newtheorem{theorem}{Theorem}[section]
\newtheorem{lemma}[theorem]{Lemma}
\newtheorem{proposition}[theorem]{Proposition}
\newtheorem{corollary}[theorem]{Corollary}
\theoremstyle{definition}
\newtheorem{definition}[theorem]{Definition}
\theoremstyle{remark}
\newtheorem{remark}[theorem]{Remark}
\title[From Powers--St\o rmer to L\"owner inertia]{From Powers--St\o rmer to L\"owner Inertia: Weighted Trace Functionals and Orthogonal Perturbations}
\author{Trung Hoa Dinh}
\thanks{ORCID: T. H. Dinh, \href{https://orcid.org/0000-0001-6303-1427}{0000-0001-6303-1427}; C. T. Le, \href{https://orcid.org/0000-0003-1323-3076}{0000-0003-1323-3076}.}
\address{Laboratory for Applied and Industrial Mathematics, Institute for Computational Science and Artificial Intelligence, Van Lang University, Ho Chi Minh City, Vietnam}
\address{Department of Mathematics and Statistics, Troy University, Troy, Alabama 36082, USA}
\email{hoa.dinhtrung@vlu.edu.vn}
\author{Cong Trinh Le}
\address{Department of Mathematics and Statistics, Quy Nhon University, Quy Nhon, Vietnam}
\email{lecongtrinh@qnu.edu.vn}
\author{Minh Toan Ho}
\address{Institute of Mathematics, Vietnam Academy of Science and Technology, Hanoi, Vietnam}
\email{hmtoan@math.ac.vn}
\subjclass[2020]{Primary 15A45, 47A63; Secondary 15A60, 47A64, 46L07}
\keywords{Powers--St\o rmer's inequality, Ricard's inequality, operator convexity, orthogonal perturbation, Mellin transform, L\"owner matrix, inertia}
\begin{document}
\begin{abstract}
We study the weighted trace functional
$\TT_{f,h}(A,B)=\Tr h^\circ(A-B)[f(A)-f(B)]$, where
$h^\circ(H)=h(H_+)-h(H_-)$. The Jordan decomposition of $A-B$ yields
exact common-roof decompositions for operator monotone and
operator convex functions, including the Powers--St\o rmer and Ricard
trace comparisons. The orthogonal cross term gives a finite rank-one
realization of the classical conditional L\"owner-matrix criterion
for operator convexity. For powers, resolvent squares yield a
coupling-dependent Mellin strip, endpoint identities, and quantitative
bounds on off-sector spectral interactions. Linearization produces
positive and inverse moment subspaces of power L\"owner matrices.
A confluent determinant calculation provides the remaining sign data
in the known inertia classification. We give explicit equality conditions
and a finite-copy application of the Powers--St\o rmer equality case.
\end{abstract}
\maketitle
\setcounter{tocdepth}{1}
\tableofcontents
\medskip

\section{Powers--St\o rmer, Ricard, and weighted trace inequalities}
\label{sec:intro}

The Powers--St\o rmer inequality, originally established in the
study of states of operator algebras \cite{PowersStormer},
compares the trace overlap of positive matrices with the trace
norm of their difference. A familiar form is
\begin{equation}\label{eq:ps}
 \frac12\Tr(X+Y-|X-Y|)
 \le \Tr X^\alpha Y^{1-\alpha},
 \qquad X,Y\ge0,\quad 0<\alpha<1.
\end{equation}
Extensions involving nonnegative operator monotone functions
were obtained by Dinh--Osaka--Ho \cite{DOH2013}; see also
Ogata \cite{Ogata2011}. The inequality has applications to
quantum state discrimination \cite{Audenaert2007}.

For nonnegative scalars, one has
\[
 (a-b)(a^{p-1}-b^{p-1})\ge |a-b|^p,
 \qquad a,b\ge0,\quad p\ge2.
\]
Ricard established the corresponding noncommutative trace
inequality \cite{Ricard2016}:
\begin{equation}\label{eq:ricard}
 \Tr|A-B|^p\le T_p(A,B),\qquad
 T_p(A,B):=\Tr(A-B)(A^{p-1}-B^{p-1}),
 \quad p\ge2.
\end{equation}
On the other hand, applying \eqref{eq:ps} to $A^p$ and $B^p$
with $\alpha=1/p$, and then interchanging $A$ and $B$, gives
\begin{equation}\label{eq:upper}
 T_p(A,B)\le\Tr|A^p-B^p|,\qquad p>1.
\end{equation}
This comparison and related inequalities for operator convex
functions were studied by Dinh--Ho--Le--Vo \cite{DHLV2019}.
Consequently, for $p\ge2$,
\begin{equation}\label{eq:chain}
 \Tr|A-B|^p\le T_p(A,B)\le\Tr|A^p-B^p|.
\end{equation}

Our starting point is a weighted extension of the trace
functional in \eqref{eq:ricard}. Let $H=A-B$ and let
$h:[0,\infty)\to[0,\infty)$ be nondecreasing with $h(0)=0$.
Writing $H=H_+-H_-$ for its Jordan decomposition, define
$h^\circ(H)=h(H_+)-h(H_-)$. For a continuous real-valued
function $f$ satisfying $f(0)=0$, we consider
\[
 \TT_{f,h}(A,B)
 =\Tr h^\circ(A-B)(f(A)-f(B))
\]
and compare it with
\[
 \BB_{f,h}(A,B)=\Tr h(|A-B|)f(|A-B|).
\]
The choice $h(t)=t$ and $f(t)=t^{p-1}$ recovers the trace
functional in Ricard's inequality. More generally, the
weight $h$ allows different spectral contributions of
$A-B$ to be treated separately.

The Jordan decomposition provides a common upper bound
for $A$ and $B$. Set
\begin{equation}\label{eq:roof}
 H=A-B,\qquad
 P=H_+,\quad Q=H_-,\quad PQ=0,\qquad
 Z=A+Q=B+P.
\end{equation}
The matrix $Z$ is the common upper bound used in the trace
inequalities of \cite{DHLV2019}. We use it to obtain exact
remainder decompositions for operator monotone and operator
convex functions. An important term in these decompositions
is the orthogonal cross increment
\begin{equation}\label{eq:cross-def}
 \XX_f(B;P,Q)
 =-\Tr Q\bigl(f(B+P)-f(B)\bigr),
 \quad B>0,\quad P,Q\ge0,\quad PQ=0.
\end{equation}
Although $P$ and $Q$ have orthogonal supports, this
increment need not vanish. Its behavior depends on how
$B$ couples the two subspaces.

The connection between this cross increment and operator
convexity is closely related to classical results on
L\"owner matrices. L\"owner characterized operator monotone
functions through divided differences \cite{Loewner1934},
and Kraus established a corresponding criterion for matrix
convexity \cite{Kraus1936}. Bendat--Sherman studied the
relationship between operator monotonicity and convexity
\cite{BendatSherman1955}, while Hansen--Pedersen developed
operator Jensen inequalities \cite{HansenPedersen1982}.
Further background is available in
\cite{Bhatia1997,Hiai2010,Simon2019}.
In particular, Bhatia--Sano and Hiai--Sano established
conditional-definiteness criteria for L\"owner matrices
associated with operator convex functions
\cite{BhatiaSano2009,BhatiaSano2010,HiaiSano2012}.
We express the corresponding quadratic forms as first
variations of finite orthogonal perturbations.

A related problem concerns the inertia of power L\"owner
matrices, which was determined by Bhatia--Friedland--Jain
\cite{BFJ2016}. Fr\'echet derivatives of matrix powers
were studied by Bhatia--Holbrook \cite{BhatiaHolbrook2000},
and Bhatia--Jain subsequently investigated the distinct
inertia problem for Kwong matrices \cite{BhatiaJain2024}.
In the present setting, resolvent representations connect
the cross increments with moment conditions on L\"owner
quadratic forms. These conditions determine definite
subspaces, while a confluent determinant calculation
recovers the remaining signs in the known inertia
classification.

The main results of this paper are exact weighted remainder
identities for operator monotone and operator convex
functions, together with equality conditions and
quantitative estimates for power functions. The orthogonal
cross increment provides a finite rank-one formulation of
the conditional L\"owner-matrix criterion for operator
convexity. For matrix powers, we derive a resolvent--Mellin
representation whose convergence strip is determined by
positive and inverse coupling depths. This representation
also yields endpoint identities, logarithmic derivatives,
and bounds on spectral interactions between orthogonal
subspaces. Linearization then connects these results with
moment filtrations of power L\"owner matrices.

\begin{figure}[htbp]
\centering
\begin{tikzpicture}[
  >=Latex,
  box/.style={
    draw,
    rounded corners,
    align=center,
    text width=8.6cm,
    minimum height=9mm,
    inner sep=3.5mm
  },
  node distance=8mm
]
\node[box] (a)
{Powers--St\o rmer and\\ Ricard's inequalities};

\node[box, below=of a] (b)
{Weighted trace functionals\\
and exact common-roof remainders};

\node[box, below=of b] (c)
{Orthogonal perturbations\\
and conditional L\"owner forms};

\node[box, below=of c] (d)
{Resolvent--Mellin representations\\
and coupling depth};

\node[box, below=of d] (e)
{Moment subspaces\\
and L\"owner inertia};

\draw[->] (a) -- (b);
\draw[->] (b) -- (c);
\draw[->] (c) -- (d);
\draw[->] (d) -- (e);
\end{tikzpicture}
\caption{Connections between trace inequalities,
weighted remainders, orthogonal perturbations,
Mellin representations, and L\"owner inertia.}
\label{fig:connections}
\end{figure}
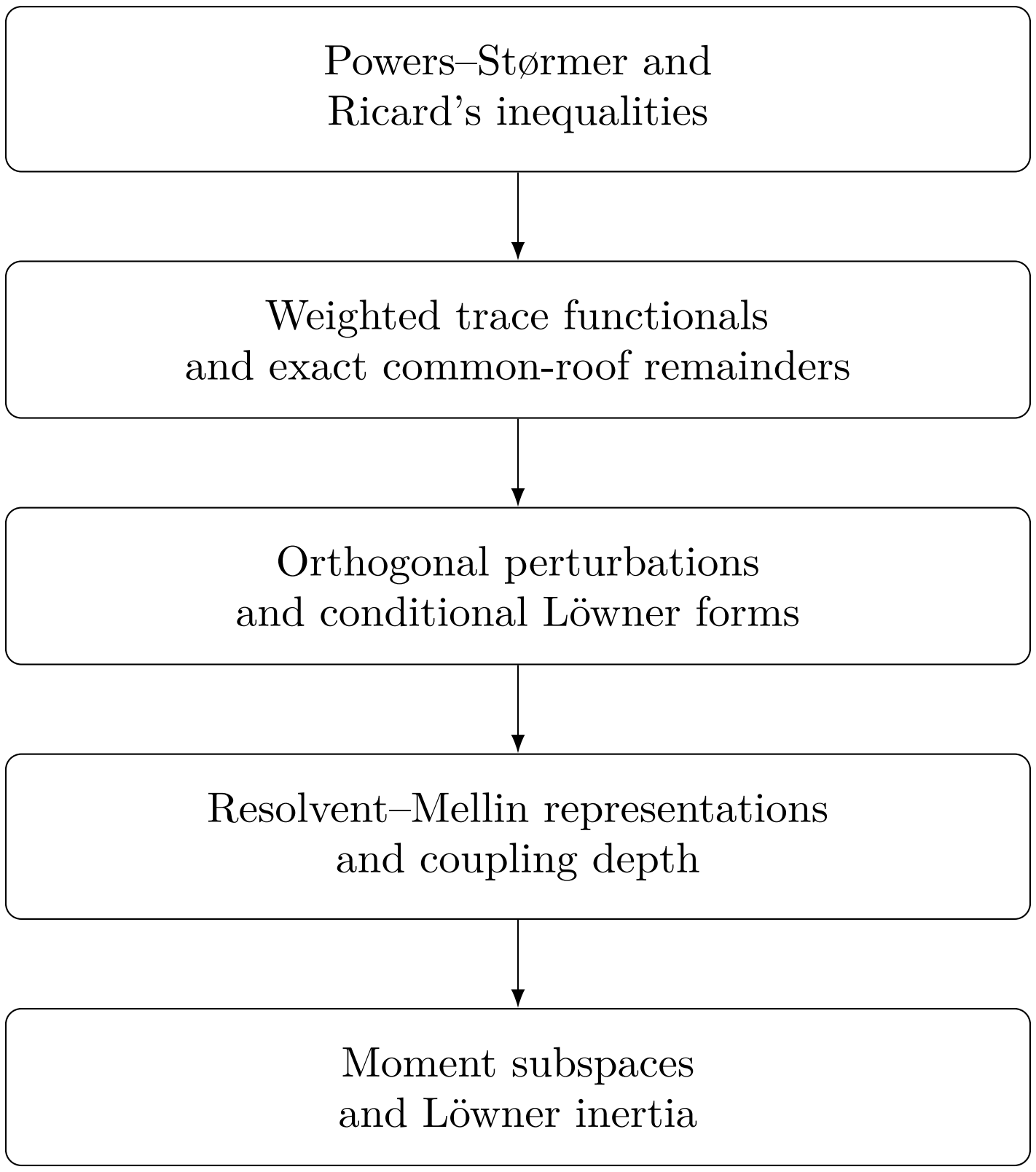

Figure~\ref{fig:connections} summarizes the connections
between the trace inequalities, orthogonal perturbations,
Mellin representations, and L\"owner inertia.

The paper is organized as follows.
Section~\ref{sec:roof} establishes the weighted trace
comparisons and exact remainder decompositions, including
equality conditions and quantitative estimates.
Section~\ref{sec:convex} relates orthogonal perturbations
to operator convexity through conditional L\"owner forms.
Section~\ref{sec:mellin} develops the resolvent--Mellin
representation and coupling-depth results.
Section~\ref{sec:loewner} studies moment subspaces and
their connection with the known inertia classification
of power L\"owner matrices.

\section{Jordan decomposition and exact trace remainders}\label{sec:roof}

\subsection*{The weighted trace functional}
For the rest of the section let $h:[0,\infty)\to[0,\infty)$ be continuous,
nondecreasing, and satisfy $h(0)=0$. If $H=H^*$, then its Jordan
parts are $H_+=(|H|+H)/2$ and $H_-=(|H|-H)/2$; they satisfy
$H_+H_-=0$. Define
\begin{equation}\label{eq:weighted-functional}
 \TT_{f,h}(A,B):=\Tr\bigl(h(H_+)-h(H_-)\bigr)
                         \bigl(f(A)-f(B)\bigr),\qquad H=A-B,
\end{equation}
and the comparison functional
\begin{equation}\label{eq:weighted-benchmark}
 \BB_{f,h}(A,B):=\Tr h(|H|)f(|H|).
\end{equation}
Here $f$ is real-valued and continuous on $[0,\infty)$ with $f(0)=0$;
we do not assume $f\ge0$ unless explicitly stated. Because $P=H_+$
and $Q=H_-$ have orthogonal supports,
\begin{equation}\label{eq:weighted-split}
\begin{split}
 \TT_{f,h}(A,B)&=\Tr h(P)[f(A)-f(B)]
                -\Tr h(Q)[f(A)-f(B)],\\
 \BB_{f,h}(A,B)&=\Tr h(P)f(P)+\Tr h(Q)f(Q).
\end{split}
\end{equation}
In particular, $h(t)=t$ gives
\[
\TT_{f,t}(A,B)=\Tr(A-B)(f(A)-f(B)),\qquad
\BB_{f,t}(A,B)=\Tr|A-B|f(|A-B|).
\]
The 2019 inequalities of Dinh--Ho--Le--Vo \cite{DHLV2019}
are recovered by setting $h(t)=t$.

For $h(t)=t^q$ and $f(t)=t^r$ with $q,r>0$ the functional becomes
\begin{equation}\label{eq:two-index-functional}
 \TT_{r,q}(A,B)=\Tr\bigl(\operatorname{sgn}(A-B)|A-B|^q\bigr)
              (A^r-B^r),
 \qquad \BB_{r,q}(A,B)=\Tr|A-B|^{q+r}.
\end{equation}
The case $q=1$, $r=p-1$ is exactly the trace quantity in the
Powers--St\o rmer--Ricard chain. Keeping $q$ free separates the power
of the spectral weight from the power in the matrix function. 
We use opposite signs for the operator monotone and operator convex
remainders:
\[
\mathscr R^{\mathrm{mon}}_{g,h}:=\BB_{g,h}-\TT_{g,h},\qquad
\mathscr R^{\mathrm{cvx}}_{f,h}:=\TT_{f,h}-\BB_{f,h}.
\]
Their positivity is governed by different classes of functions. The
common-roof decompositions contain the orthogonal cross increment
whose rank-one variation yields a L\"owner quadratic form.

\subsection*{The scalar case and spectral weights}
For commuting matrices the weighted comparison reduces to a scalar
inequality. Let $a=b+d$ with $b,d\ge0$.
For the power choice $f(t)=t^r$ and $h(t)=t^q$,
\begin{equation}\label{eq:scalar-weighted-gap}
 \mathscr T_{r,q}(a,b)-\mathscr B_{r,q}(a,b)
 =d^q\bigl((b+d)^r-b^r-d^r\bigr).
\end{equation}
When $r>1$, the difference is nonnegative; when $0<r<1$, its sign
is reversed. Indeed, for $b,d>0$,
\begin{equation}\label{eq:scalar-integral}
 (b+d)^r-b^r-d^r
 =r\int_0^b\bigl((s+d)^{r-1}-s^{r-1}\bigr)\,ds.
\end{equation}
The integral is strictly positive if $r>1$, strictly negative if
$0<r<1$, and zero when $r=1$. The transition at $r=1$ already occurs in the scalar case.
For matrices, the cross-sector terms in the remainder decomposition
can have signs governed by operator convexity, even when the total
remainder has a fixed sign.

The weight $h$ changes the relative contribution of large and
small spectral differences without changing the matrix function $f$. In the matrix setting, this distinction becomes
the separation between the function class controlling each remainder
and the spectral weight $h(P)$ or $h(Q)$ appearing in its trace.
The weighted remainder identities below separate the assumptions
on $f$ from the nonnegative spectral weight $h$.

Write
\[
 m(X,Y)=\tfrac12\Tr(X+Y-|X-Y|),\qquad
 \Delta_\alpha(X,Y)=\Tr X^\alpha Y^{1-\alpha}-m(X,Y).
\]
By Powers--St\o rmer, $\Delta_\alpha\ge0$ for $0<\alpha<1$. The following identity separates three nonnegative contributions.

\begin{proposition}\label{prop:PSroof}
For $X,Y\ge0$, let $U=(X-Y)_+$, $V=(X-Y)_-$ and $W=X+V=Y+U$. If $0<\alpha<1$, then
\begin{align}\label{eq:ps-three}
\Delta_\alpha(X,Y)
={}&\Tr(W^\alpha-Y^\alpha)Y^{1-\alpha}
+\Tr X^\alpha(W^{1-\alpha}-X^{1-\alpha})\notag\\
&+\Tr(W^\alpha-X^\alpha)(W^{1-\alpha}-Y^{1-\alpha}).
\end{align}
Each summand is nonnegative.
\end{proposition}
\begin{proof}
As $X,Y\le W$ and $0<\alpha<1$, all increments of the power functions in \eqref{eq:ps-three} are positive semidefinite. The trace of the product of two positive semidefinite matrices is nonnegative. Expanding the right-hand side leaves
$\Tr X^\alpha Y^{1-\alpha}+\Tr W-\Tr X-\Tr Y$. Finally, $m(X,Y)=\Tr X+\Tr Y-\Tr W$.
\end{proof}

The equality case is described by an orthogonal three-block decomposition.
\begin{lemma}\label{lem:strict}
If $X,Y\ge0$ and $0<r<1$, then
\[
 \Tr\bigl((X+Y)^r-X^r\bigr)X^{1-r}=0\quad\Longleftrightarrow\quad XY=0.
\]
\end{lemma}
\begin{proof}
By operator monotonicity of $t^r$, $D=(X+Y)^r-X^r\ge0$. If the trace is zero, then $D$ vanishes on $\supp X$. Write the standard resolvent formula
\begin{equation}\label{eq:strict-res}
 D=\frac{\sin\pi r}{\pi}\int_0^\infty\lambda^r
 \bigl[(X+\lambda I)^{-1}-(X+Y+\lambda I)^{-1}\bigr]d\lambda.
\end{equation}
The integrand is positive and continuous in $\lambda$; therefore its compression to $\supp X$ is zero for every $\lambda>0$. By the resolvent identity in its positive-square form, this implies
$Y^{1/2}(X+\lambda I)^{-1}(\supp X)=0$. Since $(X+\lambda I)^{-1}$ preserves and maps $\Ran(\supp X)$ onto itself, $Y(\supp X)=0$. Thus $XY=0$. Conversely, if $XY=0$, then $X$ and $Y$ have orthogonal supports and $D=Y^r$; the trace vanishes.
\end{proof}
\begin{corollary}\label{cor:ps-eq}
For one, equivalently every, $\alpha\in(0,1)$, equality holds in \eqref{eq:ps} if and only if
\begin{equation}\label{eq:splitting}
 X=C+U,\qquad Y=C+V,\qquad CU=CV=UV=0
\end{equation}
for three positive matrices $C,U,V$ with orthogonal supports.
\end{corollary}
\begin{proof}
In \eqref{eq:ps-three}, vanishing of the first two nonnegative terms and Lemma~\ref{lem:strict} yield $UY=0$ and $VX=0$. Decomposing into $\supp U$, $\supp V$, and their orthogonal complement shows that $X=U$ and $Y=0$ on the first subspace, $Y=V$ and $X=0$ on the second, and $X=Y=C$ on the third. This gives \eqref{eq:splitting}. Conversely, \eqref{eq:splitting} makes every term in \eqref{eq:ps-three} zero. The condition is independent of $\alpha$.
\end{proof}

The second roof is formed by $A^p,B^p$, not by $A,B$.
\begin{proposition}\label{prop:uppergap}
For every $p>1$,
\begin{equation}\label{eq:gap-upper}
 \Tr|A^p-B^p|-T_p(A,B)
 =\Delta_{1/p}(A^p,B^p)+\Delta_{1/p}(B^p,A^p).
\end{equation}
It is therefore an exact sum of six nonnegative roof contributions, obtained by applying Proposition~\ref{prop:PSroof} twice.
\end{proposition}
\begin{proof}
Use
$T_p=\Tr A^p+\Tr B^p-\Tr AB^{p-1}-\Tr BA^{p-1}$
and
$\Tr (A^p)^{1/p}(B^p)^{1-1/p}=\Tr AB^{p-1}$.
\end{proof}

For the lower gap let $f:[0,\infty)\to\R$ be continuous with $f(0)=0$; its behavior on the interior will be specified as needed. Define
\begin{align*}
 \EE_f(B,P)&=\Tr P\bigl(f(B+P)-f(B)-f(P)\bigr),\\
 \EE_f(A,Q)&=\Tr Q\bigl(f(A+Q)-f(A)-f(Q)\bigr),\\
 \XX_f(A;Q,P)&=-\Tr P\bigl(f(A+Q)-f(A)\bigr),\\
 \XX_f(B;P,Q)&=-\Tr Q\bigl(f(B+P)-f(B)\bigr).
\end{align*}
The notation $\EE_f$ denotes the one-sided remainder, while
$\XX_f$ denotes the orthogonal cross increment.
\begin{proposition}\label{prop:four}
With $P,Q,Z$ from \eqref{eq:roof},
\begin{align}\label{eq:four}
 &\Tr(A-B)(f(A)-f(B))-\Tr|A-B|f(|A-B|)\\
 &\qquad=\EE_f(B,P)+\EE_f(A,Q)
 +\XX_f(A;Q,P)+\XX_f(B;P,Q).\notag
\end{align}
\end{proposition}
\begin{proof}
Substitute $A-B=P-Q$ and
$f(A)-f(B)=[f(Z)-f(B)]-[f(Z)-f(A)]$.
Expanding the product produces four weighted increments. Subtracting $\Tr Pf(P)+\Tr Qf(Q)$ gives the formula, since $PQ=0$ and $f(0)=0$.
\end{proof}

\subsection*{Weighted majorization and four-term remainders}
The scalar weights in \eqref{eq:weighted-functional} are not assumed
to be operator monotone. Their monotonicity as scalar functions is enough
for the following trace comparison.

\begin{lemma}\label{lem:weak-major}
Suppose $P,X\ge0$, $a:[0,\infty)\to[0,\infty)$ is nondecreasing,
and $\phi:[0,\infty)\to[0,\infty)$ is nondecreasing. If the ordered
eigenvalues of $X$ are weakly majorized by those of $\phi(P)$, then
\begin{equation}\label{eq:major-weight}
 \Tr a(P)X\le \Tr a(P)\phi(P).
\end{equation}
\end{lemma}
\begin{proof}
List the eigenvalues of $P$ and $X$ in descending order as
$p_1\ge\cdots\ge p_n$ and $x_1\ge\cdots\ge x_n$. Von Neumann's
trace inequality gives
$\Tr a(P)X\le\sum_i a(p_i)x_i$. Since the coefficients
$a(p_i)$ decrease with $i$, summation by parts, together with
$\sum_{i\le k}x_i\le\sum_{i\le k}\phi(p_i)$, yields
$\sum_i a(p_i)x_i\le\sum_i a(p_i)\phi(p_i)$. This is
\eqref{eq:major-weight}.
\end{proof}

\begin{lemma}\label{lem:oneside-weight}
Let $g\ge0$ be operator monotone on $[0,\infty)$, $g(0)=0$.
For $B,P\ge0$ and every weight $h$ as above,
\begin{equation}\label{eq:mon-oneside}
 0\le \Tr h(P)[g(B+P)-g(B)]\le\Tr h(P)g(P).
\end{equation}
\end{lemma}
\begin{proof}
The lower bound follows from operator monotonicity. For $s>0$
put $g_s(t)=t/(t+s)$. The one-sided submajorization result
in Ricard's extension of Ando's inequalities \cite{Ricard2019}
gives $g_s(B+P)-g_s(B)\prec_w g_s(P)$.
The L\"owner representation of a nonnegative operator monotone
function vanishing at zero expresses $g$ as a nonnegative linear term
and a positive integral of the kernels $g_s$.
The Ky Fan norms and the triangle inequality preserve the resulting
weak-majorization bound after integration:
\[
 g(B+P)-g(B)\prec_w g(P).
\]
Applying Lemma~\ref{lem:weak-major} with $a=h$, $\phi=g$
gives the upper bound. The linear part contributes equality.
\end{proof}

\begin{lemma}\label{lem:convex-weight}
Suppose $f$ is operator convex on $[0,\infty)$ and $f(0)=0$.
For $B,P,Q\ge0$ with $PQ=0$, one has
\begin{align}\label{eq:weighted-cvx-defects}
 \Tr h(P)[f(B+P)-f(B)-f(P)]&\ge0,\\
 -\Tr h(Q)[f(B+P)-f(B)]&\ge0.\notag
\end{align}
\end{lemma}
\begin{proof}
First suppose $f$ has a finite right derivative at zero.
Its standard operator convex integral representation consists of
an affine term, a nonnegative quadratic term, and a positive
integral of $k_s(t)=t^2/(t+s)$ for $s>0$. The integrability factors are incorporated into the representing measure.
The affine contribution to both expressions vanishes because $h(P)$ is supported on $P$, whereas
$h(Q)P=0$.

For $f(t)=t^2$, the first quantity in \eqref{eq:weighted-cvx-defects}
is
\[
 \Tr h(P)(BP+PB)=2\Tr BPh(P)\ge0,
\]
and the second is zero since $PQ=0$. Here $Ph(P)\ge0$ and hence
$\Tr BPh(P)\ge0$.
For $k_s(t)=t-sg_s(t)$ the first expression is
\[
 s\Tr h(P)[g_s(P)-(g_s(B+P)-g_s(B))]\ge0
\]
by Lemma~\ref{lem:oneside-weight}. Its cross expression is
$s\Tr h(Q)[g_s(B+P)-g_s(B)]\ge0$, by operator
monotonicity of $g_s$. Positive integration proves both claims.

For a general continuous operator convex $f$ with $f(0)=0$, use
$f_\varepsilon(t)=f(t+\varepsilon)-f(\varepsilon)$.
Each $f_\varepsilon$ is operator convex and has a finite derivative
at zero, and the preceding argument applies. Functional calculus
and continuity as $\varepsilon\downarrow0$ give the result.
\end{proof}

\begin{theorem}\label{thm:weighted-roofs}
Let $A,B\ge0$, $H=A-B=P-Q$, $Z=A+Q=B+P$,
and let $h$ be continuous, nonnegative and nondecreasing with $h(0)=0$.
\begin{enumerate}[label=\arabic*.]
\item If $g\ge0$ is operator monotone on $[0,\infty)$ with $g(0)=0$, then
\begin{align}\label{eq:weighted-mon-roof}
\BB_{g,h}-\TT_{g,h}&=M_P+M_Q+C_P+C_Q,\\
M_P&=\Tr h(P)[g(P)-(g(Z)-g(B))],\notag\\
M_Q&=\Tr h(Q)[g(Q)-(g(Z)-g(A))],\notag\\
C_P&=\Tr h(P)[g(Z)-g(A)],\notag\\
C_Q&=\Tr h(Q)[g(Z)-g(B)].\notag
\end{align}
All four terms are nonnegative.
\item If $f$ is operator convex on $[0,\infty)$ with $f(0)=0$, then
\begin{align}\label{eq:weighted-cvx-roof}
\TT_{f,h}-\BB_{f,h}&=E_P+E_Q+X_P+X_Q,\\
E_P&=\Tr h(P)[f(Z)-f(B)-f(P)],\notag\\
E_Q&=\Tr h(Q)[f(Z)-f(A)-f(Q)],\notag\\
X_P&=-\Tr h(P)[f(Z)-f(A)],\notag\\
X_Q&=-\Tr h(Q)[f(Z)-f(B)].\notag
\end{align}
All four terms are nonnegative.
\end{enumerate}
\end{theorem}
\begin{proof}
Set $D_P=f(Z)-f(B)$ and $D_Q=f(Z)-f(A)$.
By \eqref{eq:weighted-split} and $f(A)-f(B)=D_P-D_Q$,
\[
\TT_{f,h}=\Tr h(P)D_P-\Tr h(Q)D_P
            -\Tr h(P)D_Q+\Tr h(Q)D_Q.
\]
Subtracting \eqref{eq:weighted-benchmark} yields
\eqref{eq:weighted-cvx-roof}, irrespective of signs; the
positivity of its terms is Lemma~\ref{lem:convex-weight},
applied to $(B,P,Q)$ and $(A,Q,P)$.
With $g$ in place of $f$, reverse the subtraction to get
\eqref{eq:weighted-mon-roof}. The first two terms are nonnegative by
Lemma~\ref{lem:oneside-weight}, and the last two by operator
monotonicity. In both cases, the identities are exact.
\end{proof}

\begin{corollary}\label{cor:two-param}
For all $A,B\ge0$ and $q>0$, define $\TT_{r,q}$ as in
\eqref{eq:two-index-functional}. Then
\begin{equation}\label{eq:phase-two-param}
\begin{array}{ll}
\TT_{r,q}(A,B)\le\Tr|A-B|^{q+r},&0<r\le1,\\[3pt]
\TT_{r,q}(A,B)\ge\Tr|A-B|^{q+r},&1\le r\le2.
\end{array}
\end{equation}
Both inequalities coincide at $r=1$. When $q=1$ and $r=p-1$,
this is the usual trace comparison in the operator monotone and
operator convex ranges. For $r>2$, universal termwise positivity fails, as shown below.
\end{corollary}
\begin{proof}
Take $h(t)=t^q$, $g(t)=t^r$ for $0<r\le1$, and
$f(t)=t^r$ for $1\le r\le2$ in the two parts of
Theorem~\ref{thm:weighted-roofs}. The functions are respectively
operator monotone and operator convex on $[0,\infty)$.
\end{proof}

Positivity of each one-sided or cross remainder is controlled
by the operator-function class of $f$. The weight $h$ changes
the contribution of each spectral sector of $A-B$.
There is no requirement that $h$ itself be operator monotone:
monotonicity of its scalar eigenvalue weights is sufficient for the
one-sided majorization argument.

For $f(t)=t^2$, the cross terms $X_P,X_Q$ vanish identically and
\eqref{eq:weighted-cvx-roof} reduces to the explicit identity
\begin{equation}\label{eq:weighted-quadratic}
 \TT_{t^2,h}(A,B)-\BB_{t^2,h}(A,B)
 =2\Tr BP\,h(P)+2\Tr AQ\,h(Q).
\end{equation}
The case $h(t)=t$ is the exact endpoint formula for Ricard's
remainder at $p=3$. Thus even when the cross increment is
invisible to quadratic polynomials, the one-sided overlap terms
retain information about how the spectral sectors of $A-B$
interact with the positive backgrounds $A$ and $B$.

\begin{corollary}
\label{cor:quadratic-weighted-equality}
Suppose $h(t)>0$ for every $t>0$, in addition to the standing
hypotheses on $h$. Then equality in
\eqref{eq:weighted-quadratic} holds if and only if
\begin{equation}\label{eq:weighted-equality}
 BP=0,\qquad AQ=0.
\end{equation}
Equivalently, $A=C+P$ and $B=C+Q$ for pairwise orthogonal
positive matrices $C,P,Q$. In particular, this equality criterion
is independent of the choice of such a weight $h$.
\end{corollary}
\begin{proof}
The two terms on the right of \eqref{eq:weighted-quadratic}
are nonnegative. If their sum is zero, then each vanishes.
Writing $D_P=Ph(P)\ge0$ gives
$\Tr BD_P=\|B^{1/2}D_P^{1/2}\|_2^2=0$.
Since $h$ is strictly positive on $(0,\infty)$,
$D_P$ has precisely the support of $P$, so $BP=0$.
The second term similarly yields $AQ=0$.
The common-roof decomposition then gives the orthogonal
three-block form as in \eqref{eq:splitting}.
Conversely, the three-block form makes both overlap traces zero.
\end{proof}

\subsection*{The quadratic endpoint}
For $f(t)=t^2$, the cross increment vanishes for every admissible $B,P,Q$, because
\[
 (B+P)^2-B^2=BP+PB+P^2
\]
and the trace of each summand against $h(Q)$ is zero. Nevertheless,
the \emph{total} weighted convex remainder generally does not vanish:
by \eqref{eq:weighted-quadratic}, it is twice the overlap of the
positive and negative Jordan sectors with their opposite backgrounds.
Thus at the quadratic endpoint, the one-sided remainder still records
information that is invisible to any orthogonal cross test.

For $f(t)=t^r$ with $1<r<2$, the nonnegative cross terms also
measure how the background $B$ can transmit a perturbation in the
$P$-sector into the orthogonal $Q$-sector. At $r=2$ this transmission
term disappears identically. For $r>2$ the sign need not remain
nonnegative, and the sign of a single cross term cannot be used as
a substitute for the sign of the \emph{total} Ricard remainder.
The two-dimensional calculation below exhibits a negative cross
term for the cubic power, while Ricard's inequality continues to hold.

\subsection*{A two-dimensional obstruction to cross positivity}
One can see directly why a weighted cross term can change sign
while the total trace inequality remains valid. Let
\[
B=\begin{pmatrix}a&b\\\bar b&c\end{pmatrix}>0,\quad
P=p\begin{pmatrix}1&0\\0&0\end{pmatrix},\quad
Q=q\begin{pmatrix}0&0\\0&1\end{pmatrix},\qquad p,q>0.
\]
Choose $q$ sufficiently small that $A=B+P-Q\ge0$.
For any weight $h$ with $h(q)>0$, define the cross remainder
$\CF_{f,h}(B;P,Q)$ as in \eqref{eq:weighted-cross-def}.
For a quadratic function the cross term vanishes, because
$P$ and $Q$ are orthogonal. For the cubic power,
\begin{equation}\label{eq:cubic-2by2}
 \CF_{t^3,h}(B;P,Q)=-p\,h(q)|b|^2.
\end{equation}
Indeed, the $(2,2)$ entry of $(B+P)^3-B^3$ is
$p|b|^2$: in the expansion into words in $B$ and $P$,
only the word $BPB$ survives in that entry. The sign is strictly
negative when $b\ne0$. Hence universal positivity of all the
four convex roof terms cannot extend to the cubic power.
For powers $t^r$ with $r>2$, the differential witness
\eqref{eq:witness} supplies the corresponding obstruction.
The negative cross term therefore coexists with the nonnegative
total Ricard gap.

The matrix entry $b$ is a coupling parameter: it moves vectors
between the otherwise orthogonal supports of $P$ and $Q$.
When $b=0$, the sectors are reducing and the cross term vanishes
for every function. When $b\ne0$, an odd power detects the
interaction. The same coupling is measured by the Mellin kernel in
Section~\ref{sec:mellin}.

For $2<p\le3$, the four terms in the Ricard decomposition are
nonnegative. The same resolvent calculation gives the following
quantitative estimate.
\begin{lemma}\label{lem:overlap}
For $B,P\ge0$ and $1<r\le2$, put $M=\|B+P\|$. If $M>0$, then
\begin{equation}\label{eq:overlap}
 \EE_{t^r}(B,P)
 \ge\frac{r(r^2-1)}{3}M^{r-2}\Tr BP^2.
\end{equation}
In particular $\EE_{t^r}(B,P)=0$ if and only if $BP=0$.
\end{lemma}
\begin{proof}
For $s>0$, let $g_s(t)=t/(t+s)$ and $k_s(t)=t^2/(t+s)=t-sg_s(t)$. Define
$R=(B+sI)^{-1}$,
$T=P^{1/2}RP^{1/2}$ and $T_0=s^{-1}P$.
The Woodbury identity gives
\begin{equation}\label{eq:phi}
 \EE_{k_s}(B,P)=s^2\bigl[\Tr\phi(T_0)-\Tr\phi(T)\bigr],
 \qquad \phi(x)=\frac{x^2}{1+x}.
\end{equation}
Indeed,
$\Tr P[g_s(B+P)-g_s(B)]=s\Tr\phi(T)$ and
$\Tr Pg_s(P)=s\Tr\phi(T_0)$; the linear parts cancel.
Since $0\le T\le T_0$ and $\phi'(x)=x(2+x)/(1+x)^2\ge0$ for $x\ge0$, integration along the segment from $T$ to $T_0$ shows $\EE_{k_s}\ge0$.

Set $D=T_0-T\ge0$. For $0\le X\le T_0$, scalar functional calculus gives
\[\phi'(X)\ge\frac{2X}{(1+\|T_0\|)^2}.\]
The trace derivative and $T+tD\ge T$ therefore imply
\begin{equation}\label{eq:phi2}
 \Tr\phi(T_0)-\Tr\phi(T)
 \ge\frac{2}{(1+\|T_0\|)^2}\Tr TD.
\end{equation}
Moreover,
\[
 T\ge\frac{P}{s+\|B\|},\qquad
 D\ge\frac{P^{1/2}BP^{1/2}}{s(s+\|B\|)},
\]
and positivity of trace products yields
$\Tr TD\ge \Tr BP^2/[s(s+\|B\|)^2]$. Combining with \eqref{eq:phi}--\eqref{eq:phi2},
\begin{equation}\label{eq:ks-est}
 \EE_{k_s}(B,P)\ge
 \frac{2s^3}{(s+\|P\|)^2(s+\|B\|)^2}\Tr BP^2.
\end{equation}
For $1<r<2$, write
\begin{equation}\label{eq:pow-int}
 t^r=\frac{\sin\pi(r-1)}{\pi}\int_0^\infty s^{r-2}k_s(t)\,ds.
\end{equation}
Integrating \eqref{eq:ks-est}, bounding $\|P\|,\|B\|$ by $M$, and evaluating the beta integral gives
\[
 \EE_{t^r}(B,P)\ge \frac{2\sin\pi(r-1)}\pi\Tr BP^2
 \int_0^\infty\frac{s^{r+1}}{(s+M)^4}\,ds
 =\frac{r(r^2-1)}3M^{r-2}\Tr BP^2.
\]
For $r=2$, the identity $\EE_{t^2}(B,P)=2\Tr BP^2$ gives the same bound. The coefficient is positive. Finally, $\Tr BP^2=\|B^{1/2}P\|_2^2$ vanishes exactly when $BP=0$; conversely, orthogonal supports give equality in the definition of $\EE$.
\end{proof}

For $D\ge0$, $U,V\ge0$, $UV=0$, and $s>0$, define
\begin{equation}\label{eq:G}
 \GK(D;U,V;s)=\left\|
 \bigl(I+U^{1/2}(D+sI)^{-1}U^{1/2}\bigr)^{-1/2}
 U^{1/2}(D+sI)^{-1}V^{1/2}
 \right\|_2^2.
\end{equation}
The Woodbury identity gives
\begin{equation}\label{eq:gres}
\Tr V\bigl[(D+sI)^{-1}-(D+U+sI)^{-1}\bigr]
=\GK(D;U,V;s).
\end{equation}
The cross increment is the resolvent-integral term in the quantitative
Ricard remainder below.
\begin{theorem}\label{thm:quant}
Let $A,B\ge0$, $Z=A+Q=B+P\ne0$ and $2<p<3$. Put
\[
 \gamma_p=\frac{p(p-1)(p-2)}3,\qquad
 \eta_p=\frac{\sin\pi(p-2)}\pi.
\]
For $\RR_p(A,B)=T_p(A,B)-\Tr|A-B|^p$,
\begin{align}\label{eq:quant}
\RR_p(A,B)\ge{}&\gamma_p\|Z\|^{p-3}\bigl(\Tr BP^2+\Tr AQ^2\bigr)\\
&+\eta_p\int_0^\infty s^{p-1}
 \bigl[\GK(B;P,Q;s)+\GK(A;Q,P;s)\bigr]ds.\notag
\end{align}
Furthermore,
\begin{equation}\label{eq:eq-r}
 \RR_p(A,B)=0\ \Longleftrightarrow\
 A=C+P,\ B=C+Q,\qquad CP=CQ=PQ=0
\end{equation}
for some $C,P,Q\ge0$. At $p=3$, one has the exact identity
\begin{equation}\label{eq:p3}
 \RR_3(A,B)=2\Tr BP^2+2\Tr AQ^2,
\end{equation}
with the same equality class. At $A=B=0$ all remainders vanish.
\end{theorem}
\begin{proof}
With $r=p-1\in(1,2)$, the integral \eqref{eq:pow-int} and $PQ=0$ give
\begin{equation}\label{eq:cross-mellin}
\XX_{t^r}(B;P,Q)=\eta_p\int_0^\infty
 s^r\GK(B;P,Q;s)\,ds.
\end{equation}
The analogous formula holds for $(A;Q,P)$. It is convergent: the resolvent difference is bounded by $s^{-1}I$ as $s\downarrow0$, and the squared resolvent interaction is $O(s^{-4})$ as $s\to\infty$, since $PQ=0$. Thus the two cross terms in \eqref{eq:four} are precisely the two integrals in \eqref{eq:quant}.
Apply Lemma~\ref{lem:overlap} to $\EE_{t^r}(B,P)$ and $\EE_{t^r}(A,Q)$; their ambient norms are both $\|Z\|$. Noting $r(r^2-1)/3=\gamma_p$ gives \eqref{eq:quant}.

If the remainder vanishes, both overlap traces vanish, hence $BP=AQ=0$. On $\supp P$ we then have $B=0,A=P$, and on $\supp Q$ we have $A=0,B=Q$. On their orthogonal complement $A=B=C\ge0$. This gives \eqref{eq:eq-r}. Conversely, the splitting immediately gives $\RR_p=0$. For $p=3$, $f(t)=t^2$ has zero cross increment by $PQ=0$, while each one-sided increment equals twice the corresponding overlap trace, giving \eqref{eq:p3}.
\end{proof}

\begin{corollary}\label{cor:spectral}
Fix $0<a<b$ and $2<p<3$. Let $t_1,\ldots,t_m$ be the distinct eigenvalues of $B$, with spectral projections $F_1,\ldots,F_m$. Let $K=K(a,b,p;B)$ be the $m\times m$ matrix
\begin{equation}\label{eq:K}
 K_{ij}=\int_a^b\frac{s^{p-1}}{(t_i+s)(t_j+s)}\,ds,
 \quad 1\le i,j\le m.
\end{equation}
Then $K$ is strictly positive definite, and with $\kappa_B=\lambda_{\min}(K)>0$,
\begin{equation}\label{eq:spectral}
 \RR_p(A,B)\ge\eta_p\frac{a}{a+\|Z\|}\kappa_B
 \sum_{j=1}^m\|P^{1/2}F_jQ^{1/2}\|_2^2.
\end{equation}
An analogous estimate holds with $(A,Q,P)$ in place of $(B,P,Q)$.
\end{corollary}
\begin{proof}
For $a\le s\le b$,
\[
 (I+P^{1/2}(B+sI)^{-1}P^{1/2})^{-1}
 \ge \frac{a}{a+\|Z\|}I,
\]
because $\|P\|\le\|Z\|$. The definition of $\GK$ consequently gives
\begin{align*}
 \GK(B;P,Q;s)&\ge\frac{a}{a+\|Z\|}
 \|P^{1/2}(B+sI)^{-1}Q^{1/2}\|_2^2\\
 &=\frac{a}{a+\|Z\|}\left\|
 \sum_{j=1}^m\frac{P^{1/2}F_jQ^{1/2}}{t_j+s}
 \right\|_2^2.
\end{align*}
The functions $(t_j+s)^{-1}$ are linearly independent on $[a,b]$ because the $t_j$ are distinct. Hence their weighted Gram matrix $K$ is positive definite. Integrating the last inequality against $s^{p-1}ds$ and applying the smallest-eigenvalue bound for $K$ yields \eqref{eq:spectral}, in view of \eqref{eq:quant}.
\end{proof}

\begin{remark}\label{rem:threshold}
The four summands in \eqref{eq:four}, with $f(t)=t^{p-1}$, are nonnegative for $2\le p\le3$. At $p=2$ the cross terms vanish and the one-sided terms vanish; at $p=3$ the cross terms vanish but the one-sided terms record \eqref{eq:p3}. Above $p=3$ their termwise sign cannot be maintained universally even though Ricard's \emph{total} gap remains nonnegative.

Indeed, fix $r>2$, $B=\diag(1,2)$, $u=(1,1)^T/\sqrt2$, $v=(1,-1)^T/\sqrt2$, $P=\varepsilon uu^*$, $Q=vv^*$. Then
\begin{equation}\label{eq:witness}
 \left.\frac{d}{d\varepsilon}\langle v,(B+\varepsilon uu^*)^rv\rangle\right|_{\varepsilon=0}
 =\frac14\left[r(1+2^{r-1})-2(2^r-1)\right]>0.
\end{equation}
The strict inequality follows because $t\mapsto rt^{r-1}$ is strictly convex on $[1,2]$: its endpoint average exceeds its integral average. Thus $\XX_{t^r}(B;P,Q)<0$ for small $\varepsilon>0$. The individual cross term can therefore be negative although
the global Ricard gap remains nonnegative.
\end{remark}

\begin{corollary}\label{cor:ten}
For $2\le p\le3$, the total difference between the outer terms of \eqref{eq:chain} has the exact decomposition
\begin{align*}
 \Tr|A^p-B^p|-\Tr|A-B|^p
 ={}&\Delta_{1/p}(A^p,B^p)+\Delta_{1/p}(B^p,A^p)\\
 &+\EE_{t^{p-1}}(B,P)+\EE_{t^{p-1}}(A,Q)\\
 &+\XX_{t^{p-1}}(A;Q,P)+\XX_{t^{p-1}}(B;P,Q).
\end{align*}
Each $\Delta$ splits into three nonnegative summands by \eqref{eq:ps-three}, so the total gap is a sum of ten nonnegative terms associated with the $A,B$ and $A^p,B^p$ spectral levels. Equality holds, for $p>2$, precisely for the block splitting \eqref{eq:splitting} with $X=A,Y=B$.
\end{corollary}
\begin{proof}
Add \eqref{eq:gap-upper} and \eqref{eq:four} with $f(t)=t^{p-1}$. Nonnegativity is supplied by Proposition~\ref{prop:PSroof}, Lemma~\ref{lem:overlap}, and the resolvent formula \eqref{eq:cross-mellin}, including the endpoint $p=3$. Equality for $p>2$ follows from the overlap--coupling theorem and conversely from the orthogonal splitting. At $p=2$, the lower gap vanishes identically; outer equality follows
from the upper gap.
\end{proof}

\begin{remark}
The coefficient in \eqref{eq:quant} depends on $\|Z\|$.
The Gram-matrix eigenvalue $\kappa_B$ in \eqref{eq:spectral}
can approach zero as distinct eigenvalues of $B$ coalesce.
Consequently, the estimate does not give a dimension-free lower bound
for the commutator $[A,B]$ without additional spectral assumptions.
\end{remark}

\subsection*{The two common roofs}
At the original scale, the Jordan parts
$P,Q$ describe which directions favor $A$ and which favor $B$;
the four terms in \eqref{eq:four} measure the losses incurred
when the powers of $A$ and $B$ are compared through their roof $Z$.
At the transformed scale one applies the same construction to
$A^p$ and $B^p$; the positive and negative parts of $A^p-B^p$
are usually not powers of $P,Q$, so these roofs need not coincide.
The two roof decompositions therefore distinguish the original
Jordan sectors from those of $A^p-B^p$.

If $A$ and $B$ commute, choose a common eigenbasis with eigenvalues
$a_i,b_i$. The weighted trace functional becomes
\[
 \TT_{f,h}(A,B)=\sum_i
 \operatorname{sgn}(a_i-b_i)h(|a_i-b_i|)
          [f(a_i)-f(b_i)].
\]
Both the comparisons and their remainders then reduce to scalar
statements. For noncommuting matrices, the cross terms need not vanish despite
$PQ=0$: the resolvent of $B$ can couple the supports of $P$ and $Q$.

\subsection*{Equality and finite-copy discrimination}
The equality class in the Powers--St\o rmer estimate is independent
of its exponent. In \eqref{eq:splitting}, the matrices $C,P,Q$ act
on mutually orthogonal subspaces: the common part $C$ is identical
under both hypotheses, while $P$ and $Q$ live in disjoint sectors.
This has a direct interpretation in quantum discrimination.

Let $\rho,\sigma$ be density matrices with equal prior probability
and suppose
\begin{equation}\label{eq:states-split}
\rho=C+P,\qquad \sigma=C+Q,\qquad
CP=CQ=PQ=0,\qquad c=\Tr C.
\end{equation}
For every $0<\alpha<1$, orthogonality gives
\begin{equation}\label{eq:states-overlap}
\Tr\rho^\alpha\sigma^{1-\alpha}=\Tr C=c.
\end{equation}
The Helstrom formula for $n$ identical copies is
\[
P_{e,n}^*=\frac12-\frac14
   \|\rho^{\otimes n}-\sigma^{\otimes n}\|_1.
\]
Since the tensor products share the block $C^{\otimes n}$,
and every other tensor block supported by $\rho^{\otimes n}$ is
orthogonal to every other tensor block supported by
$\sigma^{\otimes n}$, one obtains
\begin{equation}\label{eq:exact-quantum}
 \|\rho^{\otimes n}-\sigma^{\otimes n}\|_1=2(1-c^n),
 \qquad P_{e,n}^*=\tfrac12c^n.
\end{equation}
The Powers--St\o rmer overlap upper bound gives the same number,
$\tfrac12(\Tr\rho^\alpha\sigma^{1-\alpha})^n=\tfrac12c^n$.
Only the common tensor block $C^{\otimes n}$ contributes to the
optimal error probability.
Equation~\eqref{eq:exact-quantum} is the exact discrimination error
for this orthogonal three-block class; the general Chernoff theory
is developed in \cite{Helstrom,Audenaert2007}.

\section{Finite cross tests and operator convexity}\label{sec:convex}
For the finite cross test we assume $B>0$ and require $f$ only
on $(0,\infty)$. For a real $C^1$ function $f$, its divided difference is
\[
 f^{[1]}(s,t)=\begin{cases}(f(s)-f(t))/(s-t),&s\ne t,\\ f'(s),&s=t,\end{cases}
\]
and its L\"owner matrix is
\[
 L_f(t_1,\ldots,t_n)=\bigl[f^{[1]}(t_i,t_j)\bigr]_{i,j=1}^n.
\]
Recall that a Hermitian $L$ is conditionally negative definite if $c^*Lc\le0$ for every $c$ with $\sum_i c_i=0$.

\begin{theorem}\label{thm:convex}
Suppose $f\in C^1(0,\infty)$ is real-valued and
\begin{equation}\label{eq:growth}
 \liminf_{t\to\infty}\frac{f(t)}t>-\infty.
\end{equation}
The following are equivalent:
\begin{enumerate}[label=\arabic*.]
\item $f$ is operator convex on $(0,\infty)$;
\item $\XX_f(B;P,Q)\ge0$ in every dimension whenever $B>0$, $P,Q\ge0$, $PQ=0$;
\item $-\langle v,[f(B+suu^*)-f(B)]v\rangle\ge0$ for every $B>0$, $s>0$, and orthogonal unit vectors $u,v$;
\item $L_f(t_1,\ldots,t_n)$ is conditionally negative definite for every $n$ and all $t_i>0$.
\end{enumerate}
\end{theorem}
\begin{proof}
Suppose $B=\diag(t_1,\ldots,t_n)$. The Fr\'echet derivative formula gives
\[
 -\left.\frac{d}{ds}\XX_f(B;suu^*,vv^*)\right|_{s=0}
 =\langle v,Df(B)[uu^*]v\rangle=c^*L_f(t)c,
 \quad c_i=\overline{v_i}u_i.
\]
Since $u\perp v$, $\sum_i c_i=0$. Conversely, given $c\ne0$ with $\sum_i c_i=0$, set $S=\sum_i|c_i|$ and take
\[
 v_i=(|c_i|/S)^{1/2},\qquad
 u_i=e^{i\arg c_i}(|c_i|/S)^{1/2}.
\]
Then $u\perp v$ are unit vectors and $\overline{v_i}u_i=c_i/S$. Hence the infinitesimal rank-one sign is exactly conditional negative definiteness. The equivalence of conditions 1 and 4 under \eqref{eq:growth}
is the criterion of Bhatia--Sano and Hiai--Sano \cite{BhatiaSano2009,HiaiSano2012}.

If the infinitesimal sign holds at every positive base, differentiate $s\mapsto\langle v,f(B+suu^*)v\rangle$ at each $s\ge0$: its derivative is nonpositive, so condition 3 follows. Conversely, condition 3 implies its
infinitesimal limit. To obtain condition 2, decompose $P$ into rank-one increments in $\ker Q$, and telescope along successive positive base points, evaluating each increment against a spectral decomposition of $Q$. Finally, condition 2 gives condition 3 by putting $P=suu^*,Q=vv^*$.
\end{proof}

The lower-growth condition is substantive. It is not supplied by the cross test because its kernel contains the quadratic polynomials.
\begin{proposition}\label{prop:kernel}
For real $f,g\in C^1(0,\infty)$, $\XX_f=\XX_g$ for every admissible triple in every dimension if and only if $f-g$ is a polynomial of degree at most two.
\end{proposition}
\begin{proof}
For $q(t)=a+bt+ct^2$, the constant term cancels, and the remaining traces contain $\Tr QP$, $\Tr Q(BP+PB+P^2)$, all zero when $PQ=0$. Conversely, let $h=f-g$ have zero cross increments. Take $B=\diag(s,t)$ for distinct positive $s,t$, $u=(1,1)^T/\sqrt2$, $v=(1,-1)^T/\sqrt2$, and differentiate at the zero perturbation. This yields
\begin{equation}\label{eq:trapezoid}
 h'(s)+h'(t)=2\frac{h(s)-h(t)}{s-t}.
\end{equation}
For $0<s<t<w$, write $h(w)-h(s)$ as the sum of its two adjacent increments and apply \eqref{eq:trapezoid} to all three intervals. The result is
\[
 (w-s)h'(t)=(w-t)h'(s)+(t-s)h'(w).
\]
Thus $h'$ coincides on every interval with its affine interpolant and must be affine; hence $h$ is quadratic.
\end{proof}

\begin{remark}
The function $f(t)=-t^2$ has identically zero cross increment but is not operator convex. The growth condition in Theorem~\ref{thm:convex} excludes precisely this obstruction to the converse. The theorem supplies a finite interpretation of an established L\"owner-matrix criterion; the criterion itself is not a new characterization originating in this paper.
\end{remark}

For an operator convex function, the integral representation yields
a positive combination of the resolvent squares in \eqref{eq:gres}.
For powers, this kernel can be integrated explicitly in the Mellin
parameter; see Section~\ref{sec:mellin}.

\subsection*{A two-point calculation and the classical criteria}
For two spectral points, the rank-one test has the following form.
Let $B=\diag(s,t)$, where $s,t>0$ and $s\ne t$, and put
\[
 u=2^{-1/2}(1,1)^T,\qquad v=2^{-1/2}(1,-1)^T.
\]
With $P=\varepsilon uu^*$ and $Q=vv^*$, the Fr\'echet derivative
formula gives
\begin{equation}\label{eq:two-point-convex-test}
 \left.\frac{d}{d\varepsilon}
 \mathscr C_f(B;\varepsilon uu^*,vv^*)\right|_{\varepsilon=0}
 =-\frac14\left(f'(s)+f'(t)
       -2\frac{f(s)-f(t)}{s-t}\right).
\end{equation}
The quantity in parentheses is the trapezoidal-rule defect for the
derivative of $f$. A scalar convex function need not make this
defect nonnegative for all $s,t$; the matrix test is genuinely
stronger. For instance, if $f(t)=t^3$, then
\[
 f'(s)+f'(t)-2f^{[1]}(s,t)=(s-t)^2>0,
\]
so the derivative of the cross increment is negative. More precisely,
expansion of $(B+\varepsilon uu^*)^3$ gives
\begin{equation}\label{eq:cubic-cross-example}
 \mathscr C_{t^3}(B;\varepsilon uu^*,vv^*)
 =-\varepsilon\,|\langle u,Bv\rangle|^2
 =-\frac{\varepsilon}{4}(s-t)^2.
\end{equation}
Hence the finite cross test detects the failure of operator convexity
of $t^3$ already in size two. This also explains why testing only
commuting perturbations would miss the underlying phenomenon:
for a diagonal $B$ preserving both orthogonal support spaces,
all cross increments vanish regardless of the function $f$.

The historical functional-analytic route is based on L\"owner's
operator-monotonicity theorem \cite{Loewner1934}, Kraus's divided-
difference criterion for matrix convexity \cite{Kraus1936}, and the
Bendat--Sherman transformation between monotonicity and convexity
\cite{BendatSherman1955}. Hansen--Pedersen's operator Jensen theorem
\cite{HansenPedersen1982} gives a different but closely related
formulation. The conditional-definiteness results of Bhatia--Sano
\cite{BhatiaSano2009,BhatiaSano2010}, refined by Hiai--Sano
\cite{HiaiSano2012}, are the exact ingredient in
Theorem~\ref{thm:convex}. Discussions of the hierarchy of fixed-order
matrix convexity may be found in \cite{HansenTomiyama2007,Heinavaara2018},
while the integral-representation viewpoint is developed further in
\cite{FranzHiai2014,FHR2014,Simon2019}. The finite rank-one criterion in Theorem~\ref{thm:convex} is a
reformulation of the conditional L\"owner-matrix criterion of
Bhatia--Sano and Hiai--Sano.

\subsection*{Cross-term positivity versus total-gap positivity}
For the weighted roof \eqref{eq:weighted-cvx-roof} define
\begin{equation}\label{eq:weighted-cross-def}
 \CF_{f,h}(B;P,Q):=-\Tr h(Q)[f(B+P)-f(B)],\qquad PQ=0.
\end{equation}
If $Q=vv^*$ is a rank-one projection, then $h(Q)=h(1)Q$.
For a weight with $h(1)>0$, its rank-one finite cross sign is
equivalent to the sign of $\XX_f(B;P,Q)$.
The rank-one test therefore persists for every weight with
$h(1)>0$. Setting $h(t)=t$ recovers the trace functional in
\cite{DHLV2019}.

The total difference $\TT_{f,h}-\BB_{f,h}$ need not determine
the individual sign of $\CF_{f,h}$. One-sided terms can compensate for negative cross terms. Hence
Ricard's trace inequality holds for $p\ge2$, although
termwise nonnegativity of the four-roof decomposition is guaranteed
only for $2\le p\le3$. The cross-increment sign is stronger than positivity of the
complete weighted trace gap.

\begin{proposition}\label{prop:weighted-kernel}
For every admissible $B,P,Q$ with $PQ=0$, every continuous
weight $h$ with $h(0)=0$, and any polynomial
$q_2(t)=a+bt+ct^2$,
\[
 \CF_{f+q_2,h}(B;P,Q)=\CF_{f,h}(B;P,Q).
\]
If the weights include $h(t)=t$, and $f,g\in C^1(0,\infty)$
have identical cross increments in all dimensions, then $f-g$
is a polynomial of degree at most two.
\end{proposition}
\begin{proof}
For a quadratic polynomial, the increment is a linear
combination of $P$, $BP+PB$, and $P^2$. Since $h(Q)$ is
supported on $Q$ and $QP=PQ=0$, the trace of each term
against $h(Q)$ vanishes. The converse follows from
Proposition~\ref{prop:kernel}; it suffices to use the member
$h(t)=t$ of the family.
\end{proof}

The cross increment vanishes on every polynomial of degree at most
two. The cross test is insensitive to quadratic terms, whereas
the one-sided weighted remainders retain the quadratic overlap. The lower-growth hypothesis in
Theorem~\ref{thm:convex} supplies the missing orientation of
the quadratic part, while the global weighted trace functional
retains its one-sided quadratic overlap information in
\eqref{eq:weighted-quadratic}.

\subsection*{First variation of the weighted cross increment}
The rank-one derivative also applies to the weighted cross increment. Let $B=U\diag(t_i)U^*$,
$P=\varepsilon uu^*$, and $Q=\beta vv^*$ for orthogonal
unit vectors $u,v$, where $\beta>0$. Then
$h(Q)=h(\beta)vv^*$ and
\begin{equation}\label{eq:weighted-first-var}
 \left.\frac{d}{d\varepsilon}
 \CF_{f,h}(B;\varepsilon uu^*,\beta vv^*)
 \right|_{\varepsilon=0}
 =-h(\beta)\,c^*L_f(t)c,
\end{equation}
where $c_i=\overline{(U^*v)_i}(U^*u)_i$.
Thus the first variation of the weighted cross increment is the
conditional L\"owner quadratic form. Every zero-sum vector is obtained,
up to normalization, from a pair of orthogonal rank-one directions.

Nonnegativity of the total weighted remainder does not imply
nonnegativity of each cross term. In contrast, nonnegativity of
all rank-one cross increments, together with \eqref{eq:growth},
characterizes operator convexity. The first condition concerns a
two-matrix trace comparison; the second ranges over all positive
base matrices and orthogonal rank-one directions.

\section{Mellin transforms and coupling depth}\label{sec:mellin}
Fix $B>0$ and $P,Q\ge0$ with $PQ=0$. For $z\in\C$, define
\begin{equation}\label{eq:theta}
 \Theta_z(B;P,Q)=-\Tr Q\bigl((B+P)^z-B^z\bigr),
\end{equation}
and write $G(s)=\GK(B;P,Q;s)$. Because $B,B+P$ are positive definite, $z\mapsto\Theta_z$ is entire. Its positive-kernel Mellin integral has the convergence strip
determined below.

\begin{definition}
The positive and inverse coupling depths are
\begin{align}\label{eq:depth}
 d_+&=\min\{k\ge1:P^{1/2}B^kQ^{1/2}\ne0\},\\
 d_-&=\min\{k\ge1:P^{1/2}B^{-k}Q^{1/2}\ne0\},\notag
\end{align}
with value $\infty$ if there is no such integer. In finite dimensions, either both depths are infinite or both belong to $\{1,\ldots,n-1\}$.
\end{definition}
\begin{proof}[Justification of the last assertion]
By Cayley--Hamilton, $B^{-1}$ is a polynomial in $B$ and vice versa; the spans of their powers therefore agree. Consequently, all positive coupling moments vanish if and only if all inverse coupling moments vanish. If the interaction does not vanish, the first nonzero moment appears before degree $n$ by linear dependence of the powers of $B$.
\end{proof}

For $0<r<1$, functional calculus gives
\begin{equation}\label{eq:mellin-basic}
 \Theta_r(B;P,Q)=-\frac{\sin\pi r}{\pi}
 \int_0^\infty s^rG(s)\,ds.
\end{equation}
The kernel agrees with the one in \eqref{eq:cross-mellin}.

\begin{theorem}\label{thm:depth}
Assume $d_+,d_-<\infty$. The Mellin integral
\[
 M(z)=\int_0^\infty s^zG(s)\,ds
\]
converges absolutely exactly when
\begin{equation}\label{eq:strip}
 1-2d_-<\Re z<2d_++1.
\end{equation}
It is holomorphic on this strip and satisfies
\begin{equation}\label{eq:theta-mellin}
 \Theta_z(B;P,Q)=-\frac{\sin\pi z}{\pi}M(z).
\end{equation}
At the two boundary exponents,
\begin{align}\label{eq:endpoint}
 \Theta_{2d_++1}&=-\|P^{1/2}B^{d_+}Q^{1/2}\|_2^2,\\
 \Theta_{1-2d_-}&=\bigl\|(I+P^{1/2}B^{-1}P^{1/2})^{-1/2}
 P^{1/2}B^{-d_-}Q^{1/2}\bigr\|_2^2.\notag
\end{align}
\end{theorem}
\begin{proof}
The positive resolvent identity \eqref{eq:gres} gives
\[
 G(s)=\bigl\|(I+P^{1/2}(B+sI)^{-1}P^{1/2})^{-1/2}
 P^{1/2}(B+sI)^{-1}Q^{1/2}\bigr\|_2^2.
\]
At infinity, expanding $(B+sI)^{-1}=\sum_{j\ge0}(-1)^js^{-j-1}B^j$ and using the cancellation of the coupling moments of degree less than $d_+$ yields
\begin{equation}\label{eq:asymp-plus}
 G(s)=a_+s^{-2d_+-2}+O(s^{-2d_+-3}),\qquad
 a_+=\|P^{1/2}B^{d_+}Q^{1/2}\|_2^2>0.
\end{equation}
Near zero, expanding $(B+sI)^{-1}=\sum_{j\ge0}(-1)^js^jB^{-j-1}$ gives
\begin{equation}\label{eq:asymp-minus}
 G(s)=a_-s^{2d_--2}+O(s^{2d_--1}),
\end{equation}
where
\[
 a_-=\bigl\|(I+P^{1/2}B^{-1}P^{1/2})^{-1/2}
 P^{1/2}B^{-d_-}Q^{1/2}\bigr\|_2^2>0.
\]
These expansions prove both necessity and sufficiency of \eqref{eq:strip}; on compact substrips they also provide integrable majorants and hence holomorphy. Equation \eqref{eq:mellin-basic} establishes \eqref{eq:theta-mellin} for $0<r<1$, and the identity theorem extends it throughout the strip.

For the upper endpoint $b=2d_++1$, \eqref{eq:asymp-plus} gives
$M(r)=a_+/(b-r)+O(1)$ as $r\uparrow b$, while $\sin\pi r=\pi(b-r)+o(b-r)$. Hence $\Theta_b=-a_+$. For the lower endpoint $a=1-2d_-$, \eqref{eq:asymp-minus} gives $M(r)=a_-/(r-a)+O(1)$, and $\sin\pi r=-\pi(r-a)+o(r-a)$. As $\Theta_z$ is entire, these limits are precisely \eqref{eq:endpoint}.
\end{proof}

The interval in \eqref{eq:strip} is generally larger than the classical power ranges for operator convexity or concavity. For real nonintegral $r$ in that strip,
\begin{equation}\label{eq:sign}
 \sgn\Theta_r=-\sgn(\sin\pi r),
\end{equation}
because $M(r)>0$. The sign range depends on $B,P,Q$ through their coupling depths.
\begin{remark}
Let $B=\diag(1,2,3)$ and
\[
 u=(1,-2,1)^T/\sqrt6,\qquad v=(1,1,1)^T/\sqrt3,\qquad
 P=uu^*,\quad Q=vv^*.
\]
Then $u\perp v$, $u^*Bv=0$ but $u^*B^2v\ne0$, while $u^*B^{-1}v\ne0$. Thus $d_+=2$, $d_-=1$ and the maximal strip is $-1<\Re z<5$. The power $t^r$ is not operator convex for $r>2$, yet the response in this particular configuration still has the sign prescribed by \eqref{eq:sign} up to $r=5$.
\end{remark}

At the polynomial boundary $z=2$, the response vanishes because the cross increment of a quadratic polynomial is zero. Its derivative need not vanish:
\begin{equation}\label{eq:logcross}
 \Tr Q\bigl((B+P)^2\log(B+P)-B^2\log B\bigr)
 =\int_0^\infty s^2G(s)\,ds>0
\end{equation}
when the depths are finite. Thus differentiating in the exponent detects a resolvent interaction that the quadratic power itself cannot detect. Formula \eqref{eq:logcross} is one concrete sense in which a trace remainder, a Mellin transform, and the logarithmic functions of operator theory meet.

If the depths are infinite, $P^{1/2}(B+sI)^{-1}Q^{1/2}=0$ for all $s>0$, hence $G\equiv0$ and $\Theta_z\equiv0$. In fact the $B$-cyclic subspaces generated by $\Ran P$ and $\Ran Q$ are orthogonal. Consequently, equality of the cross increment is weaker than
equality in the complete Ricard remainder.

\subsection*{Integer zeros and logarithmic derivatives}
The identity \eqref{eq:theta-mellin} also determines the derivatives
of the response at its integer zeros.
\begin{proposition}
\label{prop:logarithmic}
Assume $d_+,d_-<\infty$. For every integer $k$ satisfying
$1-2d_-<k<2d_++1$,
\begin{equation}\label{eq:integer-logarithmic}
 \Theta_k(B;P,Q)=0,\qquad
 \left.\frac{d}{dz}\Theta_z(B;P,Q)\right|_{z=k}
 =(-1)^{k+1}\int_0^\infty s^kG(s)\,ds.
\end{equation}
In particular the zeros at these integer exponents are simple
whenever the coupling depths are finite. Moreover,
\begin{equation}\label{eq:logarithmic-cross-trace}
 -\Tr Q\bigl((B+P)^k\log(B+P)-B^k\log B\bigr)
 =(-1)^{k+1}\int_0^\infty s^kG(s)\,ds.
\end{equation}
\end{proposition}
\begin{proof}
Every such $k$ lies strictly inside the strip of absolute
convergence, so $M(k)$ is finite and strictly positive.
The first identity follows from $\sin\pi k=0$ in
\eqref{eq:theta-mellin}; differentiating the holomorphic identity
at $k$ gives the second, since
$-(\pi^{-1})(\pi\cos\pi k)M(k)=(-1)^{k+1}M(k)$.
Differentiation of $X^z$ with respect to $z$ for $X>0$
is $X^z\log X$, which gives \eqref{eq:logarithmic-cross-trace}.
\end{proof}
For example, $k=2$ always belongs to the interior when the
coupling is finite. The response of the quadratic power vanishes,
but its derivative is the strictly negative quantity
$-\int_0^\infty s^2G(s)\,ds$.
The logarithmic derivative detects nonzero coupling even though
the quadratic cross increment vanishes identically.

\subsection*{Rank-one response and an exact rational kernel}
For orthogonal rank-one perturbations, the resolvent kernel is
a scalar rational function. Let $u,v$ be unit vectors with
$u\perp v$, and set $P=\alpha uu^*$, $Q=\beta vv^*$,
$\alpha,\beta>0$. The Woodbury formula reduces the resolvent
interaction to
\begin{equation}\label{eq:rank-one-G}
 G(s)=\alpha\beta\,
 \frac{|\langle u,(B+sI)^{-1}v\rangle|^2}
 {1+\alpha\langle u,(B+sI)^{-1}u\rangle}.
\end{equation}
Consequently, whenever $1-2d_-<r<2d_++1$,
\begin{equation}\label{eq:rank-one-response}
\langle v,[(B+\alpha uu^*)^r-B^r]v\rangle
=\frac{\alpha\sin\pi r}{\pi}
\int_0^\infty
\frac{s^r|\langle u,(B+sI)^{-1}v\rangle|^2}
 {1+\alpha\langle u,(B+sI)^{-1}u\rangle}\,ds.
\end{equation}
The factor $\beta$ has cancelled because the displayed expression
is normalized against the unit vector $v$. Formula
\eqref{eq:rank-one-response} is a finite, nonlinear analogue of a
L\"owner quadratic form: the rational function in its numerator
will reappear when the parameter $\alpha$ tends to zero.

The equality case is geometric. Define
\[
\mathcal K_P=\operatorname{span}\{B^j\operatorname{Ran}P:j\ge0\},\quad
\mathcal K_Q=\operatorname{span}\{B^j\operatorname{Ran}Q:j\ge0\}.
\]
By the resolvent-square representation, $G(s)=0$ for all $s>0$
if and only if $P^{1/2}(B+sI)^{-1}Q^{1/2}=0$ for all $s>0$.
Expanding at infinity shows this is equivalent to
$P^{1/2}B^jQ^{1/2}=0$ for every $j\ge0$, hence to
$\mathcal K_P\perp\mathcal K_Q$. These subspaces reduce $B$;
$P$ is supported on $\mathcal K_P$ and $Q$ on $\mathcal K_Q$.
Functional calculus then yields $\Theta_z\equiv0$ for all
complex $z$. The inverse-power moments lead to the same
criterion because $B^{-1}$ is a polynomial in $B$.
Thus the vanishing cross increment is equivalent to orthogonality
of the two cyclic subspaces generated by $P$ and $Q$.

\subsection*{Coupling depth and resolvent asymptotics}
At large $s$, the first nonzero $P$--$Q$ matrix moment of $B$
occurs at degree $d_+$. At small $s$, the first nonzero moment
of $B^{-1}$ occurs at degree $d_-$. These are two distinct
filtrations of the same pair of subspaces. The endpoint identities
in \eqref{eq:endpoint} show that the first nonzero odd-power
cross response on each side records a squared Hilbert--Schmidt
coupling; the intermediate integer powers vanish.
For $d_+=d_-=1$, the interval is $-1<\Re z<3$,
which contains the whole standard operator convex power range
but is strictly larger than the operator monotone one.
If cancellations force $d_+>1$, a particular configuration can
retain a sign law for powers beyond the universal convexity
range, without contradicting the failure of operator convexity
there. This pair-dependent interval can extend beyond the universal
operator convexity range of the power function.

\subsection*{Weighted coupling energy and its equality case}
The square-kernel representation also works with a genuine
weight $h(Q)$ in place of $Q$. For $1<r<2$, define
\begin{equation}\label{eq:weighted-kernel-G}
\begin{split}
G_h(s)&:=\bigl\|(I+P^{1/2}(B+sI)^{-1}P^{1/2})^{-1/2}\\
&\hspace{37mm}\times P^{1/2}(B+sI)^{-1}h(Q)^{1/2}\bigr\|_2^2.
\end{split}
\end{equation}
The power integral and the fact that $Ph(Q)=0$ give
\begin{equation}\label{eq:weighted-mellin-C}
 \CF_{t^r,h}(B;P,Q)=\frac{\sin\pi(r-1)}\pi
            \int_0^\infty s^r G_h(s)\,ds.
\end{equation}
The integral gives the nonnegative weighted cross remainder.
The algebra is the same as the Woodbury calculation for
$G(s)$: replace $Q^{1/2}$ by $h(Q)^{1/2}$ in the
Hilbert--Schmidt square and take the weighted trace.

\begin{theorem}
\label{thm:weighted-coupling}
Assume $B>0$, $P,Q\ge0$, $PQ=0$, $1<r<2$,
and $h(t)>0$ for $t>0$. Then
\begin{equation}\label{eq:weighted-rigidity}
\CF_{t^r,h}(B;P,Q)=0
\quad\Longleftrightarrow\quad
\mathcal K_P\perp\mathcal K_Q,
\end{equation}
with the cyclic spaces defined above.
Moreover, let $0<a<b$, and let $t_1,\dots,t_m$ be the distinct
eigenvalues of $B$, with spectral projections $F_j$.
Put
\[
K_{ij}^{(r)}=\int_a^b\frac{s^r}{(t_i+s)(t_j+s)}\,ds,
\qquad\kappa_{B,r}(a,b)=\lambda_{\min}(K^{(r)})>0.
\]
Then
\begin{equation}\label{eq:weighted-stability}
\begin{split}
\CF_{t^r,h}(B;P,Q)\ge{}&
\frac{\sin\pi(r-1)}\pi\frac{a}{a+\|P\|}
\kappa_{B,r}(a,b)\\
&\times\sum_{j=1}^m\|P^{1/2}F_jh(Q)^{1/2}\|_2^2.
\end{split}
\end{equation}
\end{theorem}
\begin{proof}
The integrand in \eqref{eq:weighted-mellin-C} is nonnegative.
If its integral vanishes, continuity implies $G_h(s)=0$ for every
$s>0$. The positive factor preceding $P^{1/2}(B+sI)^{-1}$
is invertible, so
\[
P^{1/2}(B+sI)^{-1}h(Q)^{1/2}=0,\qquad s>0.
\]
Because $h(t)>0$ for positive $t$, the ranges of $h(Q)^{1/2}$
and $Q^{1/2}$ coincide. Thus the last identity is equivalent
to $\mathcal K_P\perp\mathcal K_Q$ by the cyclic-space argument
above; the converse follows by reversing those steps.

For $a\le s\le b$, the positive inverse factor in the
Woodbury square obeys
\[
(I+P^{1/2}(B+sI)^{-1}P^{1/2})^{-1}
\ge\frac{a}{a+\|P\|}I.
\]
Consequently
\[
G_h(s)\ge \frac{a}{a+\|P\|}\Bigl\|
\sum_{j=1}^m\frac{P^{1/2}F_jh(Q)^{1/2}}{t_j+s}
\Bigr\|_2^2.
\]
The rational functions $s\mapsto(t_j+s)^{-1}$ are linearly
independent on $[a,b]$ for distinct nodes, so $K^{(r)}$
is positive definite. Integrating the squared norm over $[a,b]$
against $s^rds$ and applying the Gram lower bound yields
\eqref{eq:weighted-stability}. No uniform estimate independent
of the spectral nodes is claimed.
\end{proof}

For fixed distinct spectral nodes, the weighted defect controls
the total squared Hilbert--Schmidt norm of the off-sector blocks
$P^{1/2}F_jh(Q)^{1/2}$. The constant depends on the spectral
nodes and may tend to zero as they coalesce.

\section{Linearization, moment filtrations, and L\"owner inertia}\label{sec:loewner}
L\"owner quadratic forms are first variations of finite orthogonal
responses. Fix distinct $0<t_1<\cdots<t_n$ and put
\begin{equation}\label{eq:Lr}
 L_r(t)=\bigl[\ell_r(t_i,t_j)\bigr]_{i,j=1}^n,\qquad
 \ell_r(s,t)=\begin{cases}(s^r-t^r)/(s-t),&s\ne t,\\rs^{r-1},&s=t.\end{cases}
\end{equation}
For rank-one $P=\varepsilon uu^*$, $Q=vv^*$ and $B=U\diag(t)U^*$, set
$c_i=\overline{(U^*v)_i}(U^*u)_i$, so that $\sum_i c_i=0$. Then
\begin{equation}\label{eq:linearized}
 -\left.\frac{d}{d\varepsilon}\Theta_r(B;\varepsilon uu^*,vv^*)\right|_{\varepsilon=0}
 =c^*L_r(t)c.
\end{equation}
The same vector $c$ also carries the coupling-depth conditions.

For $m,\ell\ge1$ define
\begin{align*}
 H_m^+&=\bigl\{c:\textstyle\sum_i t_i^jc_i=0,\ j=0,\ldots,m-1\bigr\},\\
 H_\ell^-&=\bigl\{c:\textstyle\sum_i t_i^{-j}c_i=0,\ j=0,\ldots,\ell-1\bigr\},\\
 F_c(s)&=\sum_{i=1}^n\frac{c_i}{t_i+s}.
\end{align*}
The positive and inverse moment conditions translate the two
coupling depths into linear constraints on $c$:
\begin{proposition}\label{prop:moment-bridge}
Let $B=U\diag(t_1,\ldots,t_n)U^*$ with distinct $t_i>0$, let $u\perp v$ be unit vectors, and set $P=\alpha uu^*$, $Q=\beta vv^*$ with $\alpha,\beta>0$. With $c_i=\overline{(U^*v)_i}(U^*u)_i$ as above,
\[
 d_+(B;P,Q)\ge m\iff c\in H_m^+,
 \qquad d_-(B;P,Q)\ge\ell\iff c\in H_\ell^-.
\]
If a depth equals $m$ or $\ell$, the corresponding next moment $\sum_i t_i^mc_i$ or $\sum_i t_i^{-\ell}c_i$ is nonzero.
\end{proposition}
\begin{proof}
Since $P^{1/2}B^jQ^{1/2}$ has rank at most one, its vanishing is equivalent to $\langle u,B^jv\rangle=0$. In the eigenbasis of $B$, the complex conjugate of this scalar is $\sum_i t_i^jc_i$. Orthogonality gives the $j=0$ moment. Applying the same argument to $B^{-j}$ proves both equivalences, and the definition of the first nonvanishing moment proves the final statement.
\end{proof}
If $c\in H_m^+\cap H_\ell^-$, expansions at infinity and zero give
\begin{align}\label{eq:F-asymp}
 F_c(s)&=(-1)^m s^{-m-1}\sum_i t_i^mc_i+O(s^{-m-2})\quad(s\to\infty),\\
 F_c(s)&=(-1)^{\ell-1}s^{\ell-1}\sum_i t_i^{-\ell}c_i+O(s^\ell)\quad(s\downarrow0).\notag
\end{align}

\begin{theorem}\label{thm:loewner-mellin}
For $c\in H_m^+\cap H_\ell^-$ and
$1-2\ell<\Re z<2m+1$,
\begin{equation}\label{eq:LM}
 c^*L_z(t)c=\frac{\sin\pi z}{\pi}
 \int_0^\infty s^z|F_c(s)|^2\,ds,
\end{equation}
where $L_z$ denotes the matrix of divided differences of $t^z$. For real nonintegral $r$ in the strip and $c\ne0$, the sign of $c^*L_r(t)c$ is the sign of $\sin\pi r$.
\end{theorem}
\begin{proof}
For $0<r<1$, the scalar identity
\[
 \ell_r(a,b)=\frac{\sin\pi r}{\pi}
 \int_0^\infty\frac{s^r}{(a+s)(b+s)}\,ds
\]
gives \eqref{eq:LM} after summing against $\overline{c_i}c_j$. The expansions \eqref{eq:F-asymp} make the integral absolutely convergent and locally uniform in the specified strip. Both sides are holomorphic there, and analytic continuation proves \eqref{eq:LM}. If $c\ne0$, distinctness of the $t_i$ implies that the rational function $F_c$ is nonzero, so its squared integral is positive for real $r$.
\end{proof}

\begin{corollary}\label{cor:LM-end}
Let $c\in H_m^+\cap H_\ell^-$. For integers $2-2\ell\le k\le2m$,
\[
 c^*L_k(t)c=0,\qquad
 c^*L_{t^k\log t}(t)c=(-1)^k\int_0^\infty s^k|F_c(s)|^2ds.
\]
At the ends of the convergence strip,
\begin{equation}\label{eq:LM-end}
 c^*L_{2m+1}(t)c=\Big|\sum_i t_i^mc_i\Big|^2,
 \qquad
 c^*L_{1-2\ell}(t)c=-\Big|\sum_i t_i^{-\ell}c_i\Big|^2.
\end{equation}
\end{corollary}
\begin{proof}
At the indicated interior integers, substitute in \eqref{eq:LM} and differentiate in the exponent; the Mellin integral converges there. At the two endpoints the leading terms in \eqref{eq:F-asymp} produce simple poles whose residues cancel the first-order zeros of $\sin\pi z$, exactly as in Theorem~\ref{thm:depth}. The limits yield \eqref{eq:LM-end}.
\end{proof}

The definite restrictions in Theorem~\ref{thm:loewner-mellin}
do not determine the full inertia: knowing the dimensions of definite
subspaces does not determine the signs on their complements. Bhatia--Friedland--Jain used positive moment spaces in their study and solved the power L\"owner inertia problem completely \cite{BFJ2016}. The inverse-moment conditions above add a symmetric analytic view at the other endpoint of the resolvent, but by themselves do not replace the determinant information.

The remaining inertia signs follow from a confluent determinant. Let $V_r(t)$ be the $2n\times2n$ confluent evaluation matrix whose two rows at $t_i$ evaluate the functions
\[
 1,t,\ldots,t^{n-1},\quad t^r,t^{r+1},\ldots,t^{r+n-1}
\]
and their derivatives, in that order. A Schur-complement computation involving the polynomial interpolation differentiation matrix gives
\begin{equation}\label{eq:det}
 \det V_r(t)=(-1)^{\binom n2}\prod_{i<j}(t_j-t_i)^2\det L_r(t).
\end{equation}

The determinant identity follows from the polynomial interpolation
differentiation matrix. Put
\[
 V=(t_i^{j-1})_{i,j=1}^n,\qquad D=\diag(t_1,\ldots,t_n),
\]
and let $V'$ have entries $(j-1)t_i^{j-2}$. First place all evaluation rows of $V_r$ before its derivative rows; the resulting block matrix is
\[
 \widetilde V_r=\begin{pmatrix}V&D^rV\\
 V'&rD^{r-1}V+D^rV'\end{pmatrix}.
\]
Writing $\mathcal D=V'V^{-1}$, the Schur complement of the first block $V$ yields
\[
 \det\widetilde V_r=(\det V)^2
 \det\bigl(rD^{r-1}+D^r\mathcal D-\mathcal D D^r\bigr).
\]
Let $\omega_i=\prod_{k\ne i}(t_i-t_k)^{-1}$ and $\Omega=\diag(\omega_i)$. Polynomial interpolation gives, for $i\ne j$,
\[
 \mathcal D_{ij}=\frac{\omega_j}{\omega_i(t_i-t_j)}.
\]
The diagonal of $D^r\mathcal D-\mathcal D D^r$ is zero, whereas for $i\ne j$ its entry is
$\omega_j\ell_r(t_i,t_j)/\omega_i$. Consequently,
\[
 rD^{r-1}+D^r\mathcal D-\mathcal D D^r
 =\Omega^{-1}L_r(t)\Omega.
\]
Thus $\det\widetilde V_r=(\det V)^2\det L_r(t)$. Returning to the interleaved row order requires $\binom n2$ transpositions, proving \eqref{eq:det}.

For nonintegral $r>0$, the $2n$ exponents are distinct, and the exponential Chebyshev system obtained by setting $t=e^x$ determines the sign of the confluent determinant \cite{KarlinStudden1966,Vavrin1994}. In fact, for increasing exponents $\beta_1<\cdots<\beta_{2n}$ the successive Wronskians of $e^{\beta_jx}$ are positive:
\[
 W_k(x)=e^{x(\beta_1+\cdots+\beta_k)}
 \prod_{1\le a<b\le k}(\beta_b-\beta_a)>0.
\]
Therefore the confluent alternant with sorted exponents is positive. Reordering the two blocks of exponents determines the sign of every leading principal minor. The ratios of successive minors then give the signs of the diagonal pivots in an $LDL^*$ factorization. This determinant identity supplies the leading-minor signs not
specified by \eqref{eq:LM}.

More explicitly, for nonintegral $r>0$, let $q=\lfloor r\rfloor$ and write $D_j(r)=\det L_r(t_1,\ldots,t_j)$, $D_0=1$. Counting the inversions of the exponent blocks gives
\begin{equation}\label{eq:minor-sign}
 \sgn D_j(r)=(-1)^{\binom j2+N_j(q)},\qquad
 N_j(q)=\begin{cases}\binom{j-q}{2},&q\le j-2,\\0,&q\ge j-1.\end{cases}
\end{equation}
Because every $D_j(r)$ is nonzero, the no-pivot factorization has pivot signs
\begin{equation}\label{eq:pivots}
 \sgn\delta_j=\sgn(D_j/D_{j-1})=
 \begin{cases}(-1)^{j-1},&j\le q+1,\\(-1)^q,&j\ge q+2.\end{cases}
\end{equation}
This closes the known nonsingular power cases, while the integer cases are handled by the polynomial factorization below.

For example, if $1<r<2$, the pivot signs are $+,-,\ldots,-$, and $L_r$ has inertia $(1,0,n-1)$; if $2<r<3$ and $n\ge3$, the signs are $+,-,+,\ldots,+$, giving $(n-1,0,1)$. At the singular integers $1\le k<n$, the identity
\[
 \frac{s^k-t^k}{s-t}=\sum_{a=0}^{k-1}s^{k-1-a}t^a
\]
factors $L_k$ through a $k\times k$ reversal matrix; the remaining $n-k$ eigenvalues are zero. The inertia formulas are those of Bhatia--Friedland--Jain
\cite{BFJ2016}; the confluent determinant calculation supplies
the signs not fixed by the moment restrictions.

\subsection*{Two-sided moment filtrations}
If the nodes $t_i$ are distinct, the Cauchy-transform
$F_c(s)=\sum_i c_i/(t_i+s)$ is a rational function with
simple poles at $-t_i$. Vanishing of its first $m$ moments at
infinity is equivalent to
\[
\sum_i c_it_i^j=0,\quad j=0,\ldots,m-1,
\]
whereas vanishing of its first $\ell-1$ Taylor coefficients
at the origin is equivalent to
\[
\sum_i c_it_i^{-j}=0,\quad j=1,\ldots,\ell-1.
\]
The shared zero-sum condition $\sum_i c_i=0$ supplies the $j=0$
constraint in both filtrations.
These conditions have different geometry. Positive moments
record cancellation under the polynomial action of $B$;
inverse moments record cancellation under its inverse.
For the orthogonal rank-one perturbation, the coefficients
are not arbitrary but take the form
$c_i=\overline{v_i}u_i$ in the eigenbasis of $B$.
The two moment filtrations are therefore the infinitesimal
shadows of the positive and inverse coupling depths of the
nonlinear perturbation.

A dimension count makes the limitation of this statement
visible. If $0\le m<n$, the positive-moment space has
codimension $m$; the same is true of the corresponding
inverse-moment space. Their intersection may have much smaller
dimension and, in exceptional positions, can be trivial.
Hence the Mellin integral determines definiteness only on
those subspaces on which its endpoint convergence hypotheses
hold. It cannot be interpreted as a direct formula for every
vector in $\mathbb C^n$ and every exponent. The confluent determinant determines the signs on the remaining
directions.

\subsection*{The known inertia classification}
The determinant calculation yields the known classification
of Bhatia--Friedland--Jain. We include the formulas because
they locate exactly what is and is not supplied by the
Mellin--moment argument. Define
$\In(L)=(\pi(L),\zeta(L),\nu(L))$, counting positive, zero,
and negative eigenvalues.

\begin{theorem}
\label{thm:BFJ-recall}
For $0<t_1<\cdots<t_n$ and $r>0$, the matrix $L_r(t)$ has:
\begin{enumerate}[label=\arabic*.]
\item $\In(L_r)=(n,0,0)$ if $0<r<1$;
\item $\In(L_r)=(m,0,n-m)$ if $2m-1<r<2m<n$;
\item $\In(L_r)=(n-m,0,m)$ if $2m<r<2m+1<n$;
\item $\In(L_r)=(\lceil n/2\rceil,0,\lfloor n/2\rfloor)$
if $r>n-1$ and $r$ is nonintegral, and also if $r\ge n$
is an integer.
\end{enumerate}
At an integer $1\le k<n$,
\begin{equation}\label{eq:integer-inertia}
\In(L_k)=\bigl(\lceil k/2\rceil,n-k,\lfloor k/2\rfloor\bigr).
\end{equation}
For negative exponents, $L_{-r}=-D^{-r}L_rD^{-r}$,
$D=\diag(t_1,\ldots,t_n)$, so positive and negative indices
are interchanged and the nullity is unchanged.
\end{theorem}
\begin{proof}[Explanation from the determinant method]
For $r>0$ not an integer, the leading minors have the signs
in \eqref{eq:minor-sign}, hence every pivot in the $LDL^*$
factorization has the sign shown in \eqref{eq:pivots}.
Counting positive and negative pivots gives the stated
nonintegral regimes. At $1\le k<n$ one factors
\[
L_k=U_kJ_kU_k^*,\qquad
 U_k=(t_i^{j-1})_{1\le i\le n,1\le j\le k},
\]
where $J_k$ is the $k\times k$ reversal matrix.
The columns of $U_k$ are independent; $J_k$ has
$\lceil k/2\rceil$ positive and $\lfloor k/2\rfloor$
negative eigenvalues. The other $n-k$ eigenvalues vanish.
If $k\ge n$, the exponent blocks in the confluent
determinant remain disjoint and the alternating pivots give
the final nonsingular regime. Finally, the diagonal congruence
for $L_{-r}$ follows immediately from the divided difference.
This is the inertia classification of \cite{BFJ2016}.
\end{proof}

The universal cross-increment sign, subject to the growth
condition, characterizes operator convexity. For an individual
power L\"owner matrix, the Mellin formula determines definiteness
on moment subspaces, and the confluent determinant determines
the remaining pivot signs.

\subsection*{Low-dimensional inertia examples}
Consider $n=3$ with $0<t_1<t_2<t_3$. In the interval
$1<r<2$, the known inertia is $(1,0,2)$. The subspace
\[
 H_1^+(t)=\{c\in\C^3:c_1+c_2+c_3=0\}
\]
has dimension two and the L\"owner--Mellin formula shows that
$c^*L_rc<0$ for every nonzero $c\in H_1^+(t)$.
Consequently, there are at least two negative directions;
the remaining one-dimensional direction must be accounted for
by information beyond that restriction. In this particular case,
positivity on the remaining direction follows from the
principal-minor calculation, giving exactly one positive eigenvalue.
For $2<r<3$, the Mellin sign changes: the same moment subspace
is positive definite. The inertia becomes $(2,0,1)$.
Thus the parity change in the sine factor is visible in the
number of definite moment directions without by itself determining
the complementary pivot.

For $n=4$ and $3<r<4$, the space $H_2^+(t)$ has dimension two;
its restriction is negative definite. The Bhatia--Friedland--Jain
inertia is $(2,0,2)$, but the negative restriction alone yields
only the inequality $\nu(L_r)\ge2$. The missing upper bound and
the existence of positive directions come from the confluent
Vandermonde determinant. The Mellin representation determines signs on the stated
moment subspaces; the determinant gives the full inertia.
For determinant methods and their applications to nearby kernels,
see \cite{BFJ2016,Vavrin1994,Vavrin1997,BhatiaJain2015,BhatiaJain2024}.

\section*{Use of AI tools.}  ChatGPT was used for reference
exploration, preliminary calculations, drafting assistance, and language
editing. The authors are responsible for verifying all mathematical
statements, references, and conclusions.

\end{document}